**Dynamics of Interest Rate Curve by Functional Auto-Regression**

V. Kargin[1] (Cornerstone Research)
A. Onatski[2] (Columbia University)

**October, 2004**

**Abstract**
The paper uses functional auto-regression to predict the dynamics of interest rate curve. It estimates the auto-regressive operator by extending methods of the reduced-rank auto-regression to the functional data. Such an estimation technique is better suited for prediction purposes as opposed to the methods based either on principal components or canonical correlations. The consistency of the estimator is proved using methods of operator theory. The estimation method is used to analyze dynamics of Eurodollar futures rates. The results suggest that future movements of interest rates are predictable at 1-year horizons.

**1. Introduction**

With the advance of the information technology and the increasing emphasis on the quantitative analysis the data sets available to researchers become more and more extensive and detailed. An increasing number of papers use functional data, where each observation is a curve as opposed to a finite-dimensional vector, for the analysis of problems from different disciplines. This paper uses functional data on the term structure of interest rates to predict the dynamics of interest rate curve.

It is a widespread opinion that interest rate dynamics can be completely described in terms of 3 factors, which are often modeled as principal components of the interest rate variation. Indeed, more than 95% percent of the variation can be decomposed into 3 factors. However, recent research (Cochrane and Piazzesi (2002)) has suggested that this opinion may have flaws.

The new evidence indicates that projecting interest rate curve to the 3 main principal components severely handicaps the ability to predict future interest rates. In particular, it appears that the best predictors of interest rate movements are among those factors that do not contribute much to the overall interest rate variation.

If not principal components, then what statistical tool is appropriate for the analysis of interest rates predictability? One possibility is to use predictive factors, an equivalent of the simultaneous linear predictions introduced in the static context by Fortier (1966). For the time series data, the method is closely related to but is better suited for the prediction purposes than the canonical analysis of Box and Tiao (1977), which extends the classical canonical correlation analysis originally developed by Hotelling in early 1930s.

The main idea of the predictive factor method is to focus on estimation of those linear functionals of the data that can contribute most to the reduction of expected predictive error. This goal necessitates a subtle balance between the search for large correlations and large variances hidden in the data.

In the finite-dimensional case, the predictive factor analysis is related to the reduced-rank auto-regression extensively studied by Reinsel (1983). A contribution of this paper is an extension of the analysis, which accounts for the functional nature of interest rate data. It parallels in many respects the extension of the classical canonical correlation analysis to the functional data performed by Leurgans, Moyeed and Silverman (1993).

---

[1] skarguine@cornerstone.com ;
Cornerstone Research, 599 Lexington Avenue floor 43, New York, NY 10022
[2] ao2027@columbia.edu;
Economics Department, Columbia University, 1007 International Affairs Building MC 3308, 420 West 118th Street , New York, NY 10027

Our main theoretical results are in Theorems 2 and 3. Theorem 2 relates the predictive factors to eigevectors of a certain generalized eigenvalue problem. Since by the Courant-Fischer theorem, the eigenvectors of the eigenvalue problem can be characterized as solutions of a minimax problem, results of Theorem 2 suggest estimating the predictive factors as solutions of a regularized minimax problem. Theorem 3 proves that with a certain choice of the regularization parameter this procedure is consistent.

As an application, we illustrate the method using ten years of data on Eurodollar futures contracts. Consistently with previous research, we find that the best predictors of future interest rates are not among the largest principal components but are hidden among the residual components.

Our empirical analysis contributes to a long-standing problem of whether the interest rates are predictable. Some research (Duffee (2002), Ang and Piazzesi (2003)) indicates that it is hard to predict better than simply assuming random walk evolution. This means that today's interest rate is the best predictor for tomorrow's interest rate, or, for that matter, for the interest rate one year from now. The subject, however, is torn with controversy. Cochrane and Piazzesi (2002) and Diebold and Lie (2002) report, on the contrary, that their methods improve over random walk prediction. We find that the controversy may be result of whether the researchers choose to restrict attention to 3 main principal components or not.

The limitation of our approach is that we do not attempt to use non-interest rate information such as the current level and innovation in inflation, GDP growth rate and other macro variables. Recently it was discovered (Ang and Piazzesi (2003)) that there is a significant correlation in the dynamic of macro variables and interest rates. We believe that after suitable modification our method can also be applied to include macro variables.

We also do not aim to derive implications of the interest rate predictability neither for the control of economy by interest rate targeting, nor for the management of financial portfolios. We believe that methods of functional data analysis may be as useful in the problems of control as in the problems of estimation. With respect to the financial applications, we want to point out that predictability of future interest rates does not necessarily imply arbitrage opportunities. Rather the relevant question is whether portfolios that correspond to the predictable combinations of interest rates generate excess returns that cannot be explained by traditional risk factors. This is a question for a separate research effort.

The rest of the paper is organized as follows. The model is described in Section 2. The estimation problem is described in section 3. The principal component method of estimation is in Section 4. The predictive factor analysis is in Section 5. The data are described in Section 6. The results of estimation of the predictive factors are in Section 7. And Section 8 concludes.

**2. Model**

Let $P_t(T)$ denote time $t$ price of a coupon-free bond with maturity at time $t+T$. Assume this function is differentiable in $T$. Then we define forward rates as

$$f_t(T) = -\frac{\partial \log P_t(T)}{\partial T}.$$

The interest rate can conceivably exist for every positive maturity, which is less than a certain bound $\overline{T}$. But in practice, only a discrete subset of these interest rates is observable at each moment of time. In addition, the maturities of observable interest rates are typically changing from time to time and an interpolation procedure is usually used to infer the curve of the interest rates. These facts mean that a natural way to model the interest rates is as an imperfectly observed function of the continuous parameter $T$.

Let us abuse the notation and use $f_t(T)$ for the forward rate curve with the subtracted mean, $f_t(T) - \bar{f}(T)$. We will model the forward interest rates as evolving according to one of the following functional auto-regressions:

(A) $$f_{t+\delta}(T) = \rho[f_t(T)] + \varepsilon_t(T)$$

or

(B) $$f_{t+\delta}(T) - f_t(T) = \rho[f_t(T) - f_{t-\delta}(T)] + \varepsilon_t(T),$$

where $f_t(T)$ is an random variable that takes values in the real Hilbert space of square-integrable functions on $[0, \bar{T}]$, $\rho$ is a Hilbert-Schmidt integral operator, and $\varepsilon_t(T)$ is a strong H-white noise in the Hilbert space. Appendix A briefly describes the formalism of Hilbert space valued random variables and explains how it relates to a more familiar language of random processes.

The assumption that coefficient operator, $\rho$, is a Hilbert-Schmidt integral operator means that $\rho$ acts as follows

$$\rho : x(T) \to \int_0^\infty \rho(S,T) x(T) dT ,$$

where $\int_0^\infty \int_0^\infty \rho(S,T)^2 \, dS dT < \infty$.

The two models are particular cases of the general model

(M) $$f_{t+\delta}(T) = \sum_{i=0}^n (D_i + \rho_i) L^i [f_t(T)] + \varepsilon_{t+\delta}(T),$$

where $L$ is the time-shift operator:
$$L : f_t(T) \to f_{t-\delta}(T),$$
and $D_i$ and $\rho_i$ are differential and integral operators. However, we will not work in this generality.

In intuitive terms, in model (A) the expectation of future values of the forward rate curve is determined by the current values. In model (B) the expectation of the future changes is determined by the past changes. We aim to estimate operator $\rho$ using a finite sample of imperfectly observed curves and predict the future curve using this estimate. In the following sections we describe several approaches to the estimation.

## 3. The Estimation Problem

To approach the estimation of $\rho$, let us first relate it to covariance and cross-covariance operators of the functional process. Consider for definiteness model (A). Let $\Gamma_{11}$ be variance operator of random curve $f_t$ and $\Gamma_{21}$ be the cross-covariance operators for curves $f_t(T)$ and $f_{t+\delta}(T)$. From the model we know that

$$\begin{aligned}\Gamma_{12}(S,T) &= E\{f_{t+\delta}(S) f_t(T)\} \\ &= E\left\{\int \rho(S,S') f_t(S') f_t(T) dS'\right\} \\ &= \int \rho(S,S') \Gamma_{11}(S',T) dS',\end{aligned}$$

or, in operator form:

$$\Gamma_{12} = \rho \Gamma_{11}.$$

Assuming that $\Gamma_{11}$ is invertible, we have $\rho = \Gamma_{12}\Gamma_{11}^{-1}$. The natural procedure would be to substitute estimates of the covariance and cross-covariance operators into this formula. As we will see shortly, substituting the empirical covariance and cross-covariance operators will not work properly, so we will need to modify this idea.

Indeed, the kernels of the empirical covariance and cross-covariance operators are

$$\hat{\Gamma}_{11}(T_1, T_2) = \frac{1}{n}\sum_{i=1}^{n} f_{i\delta}(T_1) f_{i\delta}(T_2) - \frac{1}{n}\sum_{i=1}^{n} f_{i\delta}(T_1) \frac{1}{n}\sum_{i=1}^{n} f_{i\delta}(T_2).,$$

$$\hat{\Gamma}_{12}(T_1, T_2) = \frac{1}{n-1}\sum_{i=1}^{n-1} f_{i\delta}(T_1) f_{(i+1)\delta}(T_2) - \frac{1}{n}\sum_{i=1}^{n} f_{i\delta}(T_1) \frac{1}{n}\sum_{i=1}^{n} f_{i\delta}(T_2)..$$

The empirical covariance operator $\hat{\Gamma}_{11}$ is singular and cannot be inverted. Therefore, we must use a regularization method to obtain a consistent estimate of $\rho$.

Bosq (2000) advocates regularization by projecting on the finite number of principal components of the empirical covariance operator $\hat{\Gamma}_{11}$. Another approach is to use a functional data generalization of the canonical analysis by Box and Tiao (1977). We will discuss these methods in the next section and then propose a different method that combines their advantages.

## 4. Principal Components and Canonical Correlations

The idea of the principal component method is to determine how the operator acts on those linear combinations of $f_t(T)$ that have the largest variation. Denote the span of $k_n$ eigenvectors of $\hat{\Gamma}_{11}$ associated with the largest eigenvalues, as $H_{k_n}$, and let $\pi_{k_n}$ be the orthogonal projector on this subspace. Define the regularized covariance and cross-covariance estimates as follows: $\tilde{\Gamma}_{11} = \pi_{k_n} \hat{\Gamma}_{11} \pi'_{k_n}$ and $\tilde{\Gamma}_{12} = \pi_{k_n} \hat{\Gamma}_{12} \pi'_{k_n}$. These are simply the empirical covariance and cross-covariance operators restricted to $H_{k_n}$. Then define

$$\tilde{\rho} = \pi'_{k_n} \tilde{\Gamma}_{12} \tilde{\Gamma}_{11}^{-1} \pi_{k_n}.$$

Note that $\tilde{\rho}$ is $\tilde{\Gamma}_{12}\tilde{\Gamma}_{11}^{-1}$ on $H_{k_n}$, and zero on the orthogonal complement to $H_{k_n}$. The claim is that under certain assumptions on the covariance operator, this estimator is consistent.

Let $a_1 = (\lambda_1 - \lambda_2)^{-1}$, and $a_i = \max\{(\lambda_{i-1} - \lambda_i)^{-1}, (\lambda_i - \lambda_{i+1})^{-1}\}$ for $i > 1$, where $\lambda_i$ are eigenvalues of the covariance operator $\Gamma_{11}$ ordered in the decreasing order.

**Theorem 1**. *If for some $\beta > 1$*

$$\lambda_{k_n}^{-1} \sum_{1}^{k_n} a_j = O(n^{1/4}(\log n)^{-\beta}),$$

*then $\tilde{\rho}_n$ is consistent in operator norm induced by $L^2$ norm:*

$$\|\tilde{\rho}_n - \rho\|_{L^2} \to 0 \ a.s.$$

**Remark:** The condition of the theorem requires that the eigenvalues of the covariance matrix do not approach zero too fast, and that the eigenvalues are not too close to each other.
**Proof:** This is a restatement of Theorem 8.7 in Bosq (2000).

While asymptotically consistent, the principal component estimation method may perform very badly in small samples if the best predictors of the future evolution have little to do with the largest principal components. In this case we would be better off by searching for good predictors directly without first projecting interest rates on the largest principal components.

Statisticians long ago recognized the need for a method of finding the most important relations between two random vectors. In early 1930s Harold Hotelling, who previously invented the principal components, suggested a new method of data analysis, the canonical correlation analysis (CCA). This method is described in detail in Anderson (1984) and extended to the time series data by Box and Tiao (1977).

Here is how the canonical correlation method works in finite-dimensional situations. We are given two random vectors, $x$ and $y$. The first step is to find another pair of vectors, $u$ and $v$ that maximize covariance of scalar products, $ux$ and $vy$, subject to constraints $Var(ux) = Var(vy) = 1$. The maximizing vectors are called **the first canonical variates,** $u^{(1)}$ and $v^{(1)}$.

At the next stage we find another couple of vectors, $u$ and $v$, that solve the same maximization problem but with added constraints that $ux$ and $vy$ are uncorrelated with $u^{(1)}x$ and $v^{(1)}y$, respectively. And so on.

While definition of canonical variates can be translated straightforwardly to the functional case, the estimation requires imposing additional constraints on the function smoothness. Leurgans, Moyeeed and Silverman (1993) developed a regularization method for the estimation of functional canonical variates and proved its consistency.

To see why both the principal component and the canonical correlation methods can be deficient from the point of view of prediction, consider the following stylized example. In this example we will ignore the estimation issues and simply assume that we are able to estimate well only one of the factors. It will enable us to see the problem stripped from the technical details.

Suppose we aim to predict the income of a particular category of married couples: babysitter husbands and freelance photographer wives. Let us assume that the income of a typical babysitter follows a very persistent but not very volatile AR(1) process $x_{t+1} = ax_t + \varepsilon_{t+1}$ whereas the income of a freelance photographer follows a volatile but not particularly persistent process $y_{t+1} = by_t + \eta_{t+1}$. Assume also that the loss from a prediction $(\hat{x}_{t+1}, \hat{y}_{t+1})'$ is equal to $E\{(x_{t+1} - \hat{x}_{t+1})^2 + (y_{t+1} - \hat{y}_{t+1})^2\}$ and that we can use only one factor for prediction. What factor should we use?

This example is a two-dimensional version of the functional autoregression with

$$\Gamma_{11} = \begin{pmatrix} \dfrac{\sigma_\varepsilon^2}{1-a^2} & 0 \\ 0 & \dfrac{\sigma_\eta^2}{1-b^2} \end{pmatrix},$$

$$\Gamma_{12} = \Gamma_{21} = \begin{pmatrix} \dfrac{a\sigma_\varepsilon^2}{1-a^2} & 0 \\ 0 & \dfrac{b\sigma_\eta^2}{1-b^2} \end{pmatrix},$$

$$\rho = \Gamma_{12}\Gamma_{11}^{-1} = \begin{pmatrix} a & 0 \\ 0 & b \end{pmatrix},$$

where $\sigma_\varepsilon^2$ and $\sigma_\eta^2$ are the variances of the mutually independent innovations $\varepsilon_t$ and $\eta_t$ respectively. Our assumption about the process is that $a > b$ but $\dfrac{\sigma_\varepsilon^2}{1-a^2} < \dfrac{\sigma_\eta^2}{1-b^2}$.

Let us first use the principal component method to forecast $x$ and $y$. The eigenvector of $\Gamma_{11}$ that corresponds to the largest eigenvalue is $e_2 = (0 \ \ 1)'$. The corresponding approximation of $\Gamma_{11}^{-1}$ is $\dfrac{1-b^2}{\sigma_\eta^2} e_2 e_2'$ and therefore

$$\begin{pmatrix} \hat{x}_{t+1} \\ \hat{y}_{t+1} \end{pmatrix} = \Gamma_{12} \frac{1-b^2}{\sigma_\eta^2} e_2 e_2' \begin{pmatrix} x_t \\ y_t \end{pmatrix} = \begin{pmatrix} 0 \\ by_t \end{pmatrix},$$

so that the loss under the principal component method is

$$L_{PC} = \frac{\sigma_\varepsilon^2}{1-a^2} + \sigma_\eta^2.$$

Now let us use the canonical correlation strategy. The first canonical variates are both equal to $e_1 = (1 \ \ 0)'$. The corresponding approximation of $\Gamma_{11}^{-1}$ is $\dfrac{1-a^2}{\sigma_\varepsilon^2} e_1 e_1'$ and

$$\begin{pmatrix} \hat{x}_{t+1} \\ \hat{y}_{t+1} \end{pmatrix} = \Gamma_{12} \frac{1-a^2}{\sigma_\varepsilon^2} e_1 e_1' \begin{pmatrix} x_t \\ y_t \end{pmatrix} = \begin{pmatrix} ax_t \\ 0 \end{pmatrix},$$

so that the loss under the canonical correlation method is

$$L_{CC} = \sigma_\varepsilon^2 + \frac{\sigma_\eta^2}{1-b^2}.$$

Comparing it to $L_{PC}$, we have

$$L_{CC} \leq L_{PC} \quad \Leftrightarrow \quad \frac{b^2 \sigma_\eta^2}{1-b^2} \leq \frac{a^2 \sigma_\varepsilon^2}{1-a^2},$$

which may or may not be true. For example, even though by assumption $b^2 < a^2$, we have also assumed that $\sigma_\eta^2/(1-b^2) > \sigma_\varepsilon^2/(1-a^2)$ and if the volatility of innovations to the income of freelance photographers is large enough, we may prefer to use the principal components method.

This stylized example shows that we should search for a method that would combine advantages of principal component and canonical correlation methods. The next section develops such a method by minimizing the expectation of the squared error of prediction.

**5. Predictive Factors**

To start with, note that both the principal component and canonical correlation methods are particular ways to approximate a full-ranked $\rho$ by a reduced-rank operator. In general, a rank $k$ approximation to $\rho$ has form

$$\rho \approx AB',$$

where $A: R^k \to L^2$ and $B': L^2 \to R^k$.

We would like to find $A$ and $B'$ that minimize the mean squared error of the prediction

(RR) $\qquad\qquad\qquad E\|f_{t+1} - AB' f_t\|^2 \to \min,$

subject to normalizing constraints $B'\Gamma_{11} B = I_k$ and $A'A$ is diagonal with decreasing elements on the diagonal. Fortier (1966) considers such problem in the static context, when predictors are not the lagged values of the forecasted series, and calls the corresponding variables $B'f_t$ simultaneous linear predictions. In what follows, we will call $B'f_t$ **predictive factors** and A **predictive factor loadings**. The operator $\rho$ can be estimated by estimating $A$ and $B'$ separately:

$$\hat{\rho} = \hat{A}\hat{B}'.$$

As long as $A$ and $B'$ are estimated consistently, we can expect $\hat{\rho}$ to converge to $\rho$ in probability as the rank parameter $k$ tends to infinity.

Similarly to Hotelling's principal components and canonical variates, the predictive factors can be defined recursively. The first predictive factor[3] $b_1'f_t$ and the first predictive factor loading $a_1$ correspond to solution of (RR) for $k=1$. The second predictive factor and factor loading are defined as solving the same problem subject to an additional constraint that $b_2'$ must be orthogonal to $b_1'$ in the metric $\Gamma_{11}$, that is $b_2'\Gamma_{11}b_1 = 0$. And so on for the third, fourth, etc. factors and factor loadings.

Appendix B proves the following:

**Theorem 2**

*i) There exists a vector $B$ solving (RR). It consists of the eigenvectors of the operator pencil $\Gamma_{21}\Gamma_{12} - \lambda\Gamma_{11}$ arranged in the order of declining eigenvalues.*

*ii) Vector $A$ solving (RR) is equal to $\Gamma_{12}B$.*

*iii) The $i^{th}$ component of $B$, $b_i$, can be found recursively by maximizing $b_i'\Gamma_{21}\Gamma_{12}b_i$ subject to $b_i'\Gamma_{11}b_j = 0, j < i$ and $b_i'\Gamma_{11}b_i = 1$.*

*iv) The maximum for the problem formulated in iii) exists, is equal to the i-th eigenvalue of the operator pencil $\Gamma_{21}\Gamma_{12} - \lambda\Gamma_{11}$, and can be interpreted as the reduction in the mean squared error of forecasting due to the i-th predictive factor.*

*v) Assuming that the first $k$ eigenvalues of $\Gamma_{21}\Gamma_{12} - \lambda\Gamma_{11}$ are distinct positive numbers and that $\Gamma_{11}$ is positive definite, the solution to (RR) is unique (up to a simultaneous change in sign of $B$ and $A$.*

Theorem 2 relates the problem of finding optimal predictive factors to a generalized eigenvalue problem. Its significance is twofold. First, it relates the problem of optimal prediction to a well studied area of generalized eigenvalue problems. Second, it suggests a method for estimation of the optimal predictive factors that proceeds by solving a regularized version of the generalized eigenvalue problem.

To build intuition and explore the relationship between predictive factors on the one hand and canonical correlations and principal components on the other hand, let us return to the photographer-babysitter example. In that example, the first predictive factor is the eigenvector corresponding to the largest eigenvalue of the matrix pencil:

---

[3] In what follows, we will denote scalar products like $\int_0^{\bar{T}} x(T)y(T)dT$ as $x'y$ and functionals transforming $y$ into $x'y$ as $x'$.

$$\begin{pmatrix} \dfrac{a^2\sigma_\varepsilon^4}{(1-a^2)^2} & 0 \\ 0 & \dfrac{b^2\sigma_\eta^4}{(1-b^2)^2} \end{pmatrix} - \lambda \begin{pmatrix} \dfrac{\sigma_\varepsilon^2}{1-a^2} & 0 \\ 0 & \dfrac{\sigma_\eta^2}{1-b^2} \end{pmatrix}.$$

Therefore, if $\dfrac{b^2\sigma_\eta^2}{1-b^2} \leq \dfrac{a^2\sigma_\varepsilon^2}{1-a^2}$, the first predictive factor is equal to

$$\left(\sqrt{\dfrac{1-a^2}{\sigma_\varepsilon^2}},\ 0\right)\begin{pmatrix} x_t \\ y_t \end{pmatrix} = \sqrt{\dfrac{1-a^2}{\sigma_\varepsilon^2}}\, x_t,$$

and the first predictive factor loading is equal to $\left(a\sqrt{\sigma_\varepsilon^2/(1-a^2)},\ 0\right)$. And hence, our prediction

$$\begin{pmatrix} \hat{x}_{t+1} \\ \hat{y}_{t+1} \end{pmatrix} = \begin{pmatrix} a\sqrt{\dfrac{\sigma_\varepsilon^2}{1-a^2}} \\ 0 \end{pmatrix} \sqrt{\dfrac{1-a^2}{\sigma_\varepsilon^2}}\, x_t = \begin{pmatrix} ax_t \\ 0 \end{pmatrix}$$

coincides with the prediction under the canonical correlation analysis.

On the other hand, if $\dfrac{b^2\sigma_\eta^2}{1-b^2} \geq \dfrac{a^2\sigma_\varepsilon^2}{1-a^2}$, then the first predictive factor is equal to $\sqrt{(1-b^2)/\sigma_\eta^2}\, y_t$, and the first predictive factor loading is equal to $\left(0\ \ b\sqrt{\sigma_\eta^2/(1-b^2)}\right)$. And therefore our prediction

$$\begin{pmatrix} \hat{x}_{t+1} \\ \hat{y}_{t+1} \end{pmatrix} = \begin{pmatrix} 0 \\ b\sqrt{\dfrac{\sigma_\eta^2}{1-b^2}} \end{pmatrix} \sqrt{\dfrac{1-b^2}{\sigma_\eta^2}}\, y_t = \begin{pmatrix} 0 \\ by_t \end{pmatrix}$$

coincides with the prediction under the principal component analysis.

Hence in the above example, the prediction using the first predictive factor coincides with one using the canonical correlation analysis if the loss under the canonical correlation analysis is smaller than the loss under the principal components analysis $L_{CC} < L_{PC}$; it coincides with the prediction using the principal components analysis when the principal components loss is smaller than the loss under the canonical correlation analysis.

### 6. Consistent Estimation of Predictive Factors

It is natural to estimate operator $\rho$ by using the empirical counterpart of the reduced-rank approximation $\rho \approx AB'$. To obtain $\hat{A}$ and $\hat{B}$, one could compute the eigenvectors of $\hat{\Gamma}_{21}\hat{\Gamma}_{12} - \lambda\hat{\Gamma}_{11}$ and use theorem 2. Unfortunately, similarly to the situation with the canonical variates, such a method of estimation would be inconsistent and the corresponding estimators meaningless. It is because the predictive factors are designed to extract those linear combinations of the data that have small variance relative to their covariance with the next period's data. Linear combinations with small variance are poorly estimated and a seemingly strong covariance (in relative terms) with the next period's data may easily be an artifact of the sample.

Leurgans, Moyeed and Silverman (1993) deal with the problem for the canonical correlation analysis by introducing a penalty for roughness of the estimated canonical covariates. We use the same idea to obtain a consistent estimate of the predictive factors. To formulate the consistency theorem, we need to introduce several assumptions and notations.

**Assumption 1** $\Gamma_{11}$ is positive definite.

**Assumption 2** *Model (A) defines a weakly stationary Hilbertian process $\{f_t\}$, such that* $E\|f_t\|^4 < \infty$.

**Assumption 3** *All eigenvalues of $\Gamma_{21}\Gamma_{12} - \lambda\Gamma_{11}$ are distinct.*

Assumption 1 ensures that the data have sufficient variation to allow the consistent estimation of the predictive factors. Assumption 2 is essentially a restriction on possible AR operators $\rho$. Indeed, as Theorem 3.1 of Bosq (2000) shows, a sufficient condition for $\{f_t\}$ to be stationary is that there exists an integer $j_0 \geq 1$ such that the norm of $\rho^{j_0}$ considered as an operator in the Hilbert space is less than 1. Assumption 3 rules out the difficult situation when several predictive factors correspond to the same eigenvalue.

Let us denote the $j$-th eigenvalue and eigenvector of the operator pencils $\Gamma_{21}\Gamma_{12} - \lambda\Gamma_{11}$, $\Gamma_{21}\Gamma_{12} - \lambda(\Gamma_{11} + \alpha I)$ and $\hat{\Gamma}_{21}\hat{\Gamma}_{12} - \lambda(\hat{\Gamma}_{11} + \alpha I)$ as $\lambda_j, \lambda_{\alpha j}, \hat{\lambda}_{\alpha j}$ and $b_j, b_{\alpha j}, \hat{b}_{\alpha j}$ respectively. Here $\alpha > 0$ is a regularization parameter. We assume that the eigenvectors are normalized so that $b_j'\Gamma_{11}b_i = \delta_{ji}, b_{\alpha j}'(\Gamma_{11} + \alpha I)b_{\alpha j} = \delta_{ji}, \hat{b}_{\alpha j}'(\hat{\Gamma}_{11} + \alpha I)\hat{b}_{\alpha j} = \delta_{ji}$.

Appendix C proves the following:

**Theorem 3**

*Suppose Assumptions 1, 2, 3 hold, $\alpha \to 0$ and $(n/\log n)^{1/2}\alpha \to \infty$ as $n \to \infty$. Then for $j = 1, 2, ...$:*

    i)    $\hat{\lambda}_{\alpha j}$ *converges in probability to* $\lambda_j$

    ii)   $(\hat{b}_{\alpha j} - b_j)\Gamma_{11}(\hat{b}_{\alpha j} - b_j) \to 0$ *in probability as* $n \to \infty$.

**Remark:** Of course, what can be consistently estimated is not the eigenvector itself, but the subspace generated by this eigenvector. For this reason, statement ii) in the above theorem holds for a particular choice in the sign of the eigenvectors $\hat{b}_{\alpha j}$ and $b_j$.

Note that condition ii) implies convergence of estimates of the predictive factors in the following sense. Suppose that we estimate a predictive factor, $b_j'f_t$, where $f_t$ is chosen at random from its unconditional distribution, by $\hat{b}_{\alpha j}'f_t$. Conditionally on our estimate $\hat{b}_{\alpha j}$, a probability that the difference between the factor and its estimate is greater by absolute value than $\varepsilon$ can be bounded as follows:

$$\Pr\{(\hat{b}_{\alpha j} - b_j)f_t > \varepsilon \mid \hat{b}_{\alpha j}\} \leq \varepsilon^{-2}Var[(\hat{b}_{\alpha j} - b_j)f_t \mid \hat{b}_{\alpha j}]$$
$$= \varepsilon^{-2}(\hat{b}_{\alpha j} - b_j)\Gamma_{11}(\hat{b}_{\alpha j} - b_j).$$

According to statement ii) of theorem 3, this bound tends in probability to zero as $n \to \infty$.

Condition ii) also implies convergence in probability of our estimates of the predictive factor loadings, $\hat{a}_{\alpha j}$. Indeed, we have:

$$\hat{a}_{\alpha j} - a_j = \hat{\Gamma}_{12}\hat{b}_{\alpha j} - \Gamma_{12}b_j$$
$$= (\hat{\Gamma}_{12} - \Gamma_{12})\hat{b}_{\alpha j} + \Gamma_{12}(\hat{b}_{\alpha j} - b_j)$$

Lemma 1 from appendix C implies that the first term in the above expression tends in probability to 0. For the second term we have:

$$\Gamma_{12}(\hat{b}_{\alpha j} - b_j) = \rho \Gamma_{11}(\hat{b}_{\alpha j} - b_j)$$
$$\leq \sqrt{(\hat{b}_{\alpha j} - b_j)\Gamma_{11}(\hat{b}_{\alpha j} - b_j)} \|\rho \Gamma_{11}^{1/2}\|,$$

which tends to zero in probability according to statement ii) of Theorem 3.

In sum, Theorem 3 essentially says that both predictive factors and predictive factor loadings can be consistently estimated by maximizing a regularized Rayleigh criterion.

### 7. Data

We use daily settlement data on the Eurodollar that we obtained from Commodity Research Bureau. The Eurodollar futures are traded on Chicago Mercantile Exchange. Each contract is an obligation to deliver a 3-month deposit of $1,000,000 in a bank account outside of the United States at a specified time. The available contracts has delivery dates that starts from several first months after the current date and then go each quarter up to 10 years into the future.

The available data start in 1982, however, we use only the data starting in 1994 when the trading on 10-year contract appeared. We interpolated available data points by cubic splines to obtain smooth forward rate curves. We restricted the curve to points that are 30 days from each other to speed up the estimation. We also removed datapoints with less than 90 or more than 3480 days to expirations. That left us with 114 points per curve and 2507 valid dates.

The main difference of futures contract from the forward contract is that it settled during the entire life of the contract, while in the forward contract the payment is made only at the settlement date. This difference and variability of short-term interest rates make the values of the forward and futures contracts different. While the difference is small for short maturities, it can be significant for long maturities. There exists methods to adjust for this difference but for our illustrative purposes we will simply ignore it.

The rate on the forward contract is approximately the forward rate that we defined above. Indeed, the buyer of the contract expects to have a negative cash flow (the price of the forward contract) on the settlement date and a positive cash flow ($1,000,000) 3 months after the settlement date. He has the following alternative investment: he buys a discount bond that will pay $1,000,000 three months after the settlement date. This costs $1,000,000 $\times P_t(T + \delta)$, where $\delta$ denotes 3 months. He complements this by selling a discount bond that matures on the settlement day. If his overall investment is zero, then he is sure that on the settlement day he can make a payment of at least $1,000,000 $\times P_t(T + \delta)/P_t(T)$. By arbitrage considerations we see that the price of the forward contract is $1,000,000 times one minus the forward rate.

The next chart illustrates the evolution of Eurodollar forward rate curves.

**Figure 1.** Forward Curve Evolution

Note: The forward curves correspond to Eurodollar rates from January 1994 to December 2003. The calendar time on the left axis and the time to maturity is on the right axis. Both are measured in working days.

### 8. One year ahead prediction of the term structure

In this section, we use the predictive factor method described in section 5 to estimate our functional autoregression models (A) and (B). For the estimation and pseudo out of sample forecasts, we use the subsample of our data starting from January 1, 1994 and ending at February 28, 2001, the onset of the last recession. We find that many of our results are very sensitive to whether we include the recession period in the data sample or not. In particular, our estimates of the autoregressive operator seem to be unstable between the period of normal growth and the recession.

The goal of our estimation is a one-year-ahead prediction of the term structure of interest rates. Instead of estimating the autoregressive operator $\rho$ at daily frequency and then taking the estimate into power 252 (business days in a year) to produce our forecasts, we estimate $\rho^{252}$ directly by considering an autoregression of the term structure on its one-year lag. We estimate this autoregression using daily data. Such a method of estimation will be more robust to possible errors-in-measurement problems in our data because the bias induced by such an error will not be amplified 250 times.

Figure 2 shows our estimates of the weights of the first predictive factor and the corresponding factor loadings ($b_1$ and $a_1$ respectively, in terminology of section 5) for $\alpha = 0$, $\alpha = 0.1$, and $\alpha = 1$. In this paper we do not study how the optimal choice of the regularization parameter $\alpha$ should be made. A systematic analysis of the optimal choice is left for future work.

The horizontal axes on figure 2 correspond to maturity measured in number of 3-month periods, the longest maturity being about 30 years. As we mentioned before, without the regularization (the case $\alpha = 0$) the estimate of the predictive factor is not consistent. The estimated factor has no sense, which is clearly confirmed by the upper left graph of figure 2. As $\alpha$ grows, the weights of the factor become smoother. For $\alpha = 0.1$, the factor weights are negative for maturities less than about 1.5 year, positive for very long maturities (more than 25 years) and wiggling around zero for other maturities. For $\alpha = 1$, the factor weights look more like a linear function with positive slope.

**Figure 2**. First predictive factor weights and loadings for $\alpha = 0$, $\alpha = 0.1$, and $\alpha = 1$.

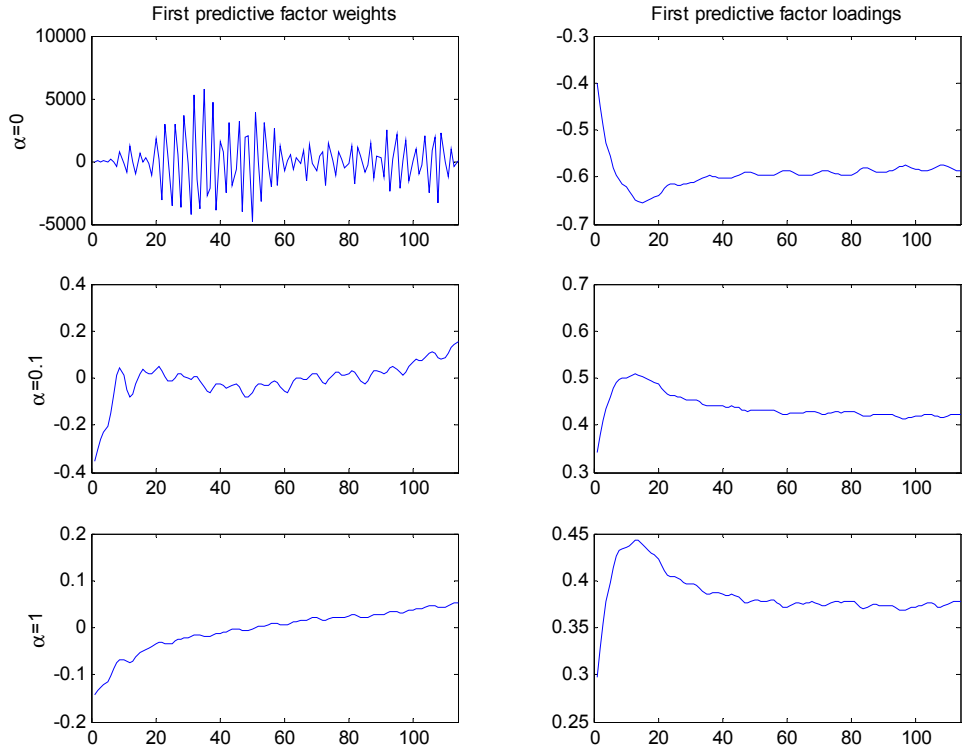

Below, we will focus on the case $\alpha = 0.1$ because we found that its pseudo out of sample forecast performance (to be described shortly) is better than that for $\alpha = 1$. Table 1 shows the first 5 eigenvalues of the operator pencil $\hat{\Gamma}_{21}\hat{\Gamma}_{12} - \lambda(\hat{\Gamma}_{11} + \alpha I)$.

**Table 1.** Eigenvalues of $\hat{\Gamma}_{21}\hat{\Gamma}_{12} - \lambda(\hat{\Gamma}_{11} + 0.1I)$.

| Eigenvalue | $\mu_1$ | $\mu_2$ | $\mu_3$ | $\mu_4$ | $\mu_5$ |
|---|---|---|---|---|---|
| | 22.03 | 0.42 | 0.04 | 0.01 | 0.00 |

According to theorem 2, the eigenvalues can be interpreted as reductions in the mean square error of forecasting due to the corresponding predictive factors. We see that the error reduction due to the first predictive factor is much larger than the reductions corresponding to the other factors. We have decided to use 3 factors to estimate the autoregressive operator mainly because of the tradition in the literature. Note, however, that the predictive factors are not designed to explain the variation in the data and hence the above reference to the tradition is not substantive. In addition, our restricting attention to 3 factors does not mean that we really believe that the rank of the autoregressive operator is equal to 3. Instead, we simply think that from the practical point of view, considering more factors in our estimation procedure would not improve the predictive value of our estimates much.

Figure 3 shows our estimates of the weights and loadings of the first three predictive factors.

**Figure 3**. The weights and loadings of the first 3 predictive factors, $\alpha = 0.1$

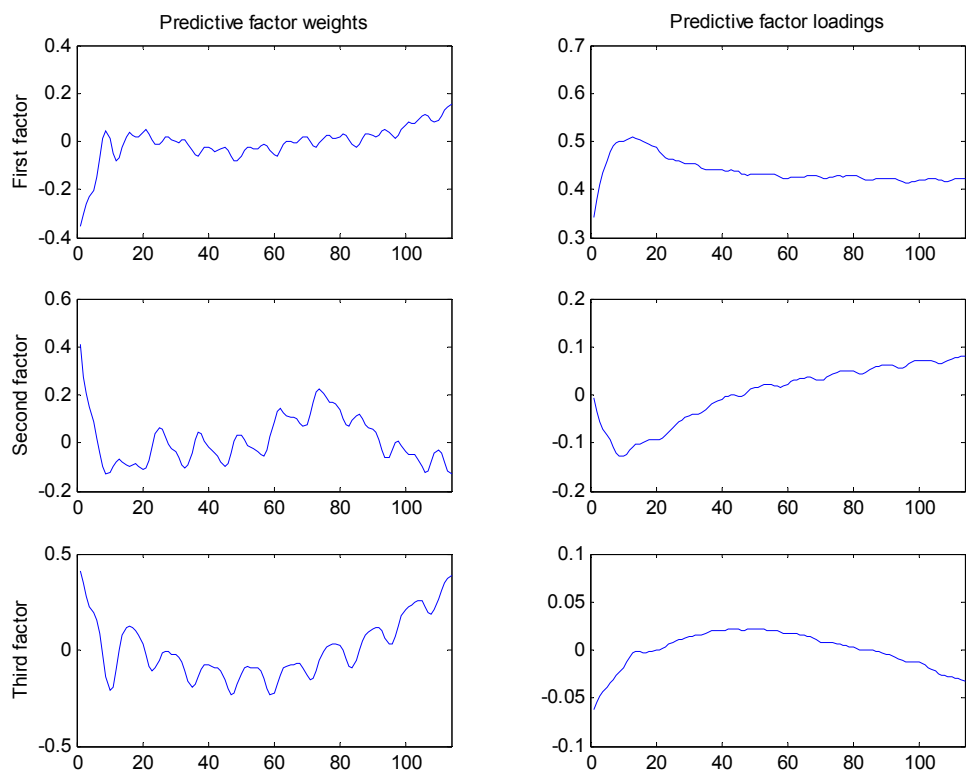

According to the estimates, an unexpected one percentage point increase in the six-month forward rate (our first observation) leads to a quarter percentage point decrease in two- year forward rates and to smaller decreases, but above 0.1 percentage points, for other maturities.

To assess the predictive performance of our estimate of model (A) based on 3 predictive factors, we run the following experiment. We first separate our sample 01/01/94:02/28/01 into two parts of equal sizes: a subsample 01/01/94:07/25/97 and a subsample 07/28/97:02/28/01. Then, we estimate our functional autoregression based on the first subsample and forecast the term structure one year ahead. The next step is to extend the first subsample to include one more day, re-estimate the functional autoregression, and forecast the term structure one year ahead and so on until we add the day one year before the end of our second subsample. After that, our forecasting would correspond to the term structures beyond the second subsample, we would not be able to compare the forecast with the actual term structure, and therefore we stop the exercise.

Our measure of the predictive performance is the root mean squared error based on the difference between actual term structure and the forecasted one. This measure will be different for different maturities. Therefore, we report a whole curve of the root mean squared errors.

We compare predictive performance of our method with 4 different methods. The first one is the same functional autoregression but estimated based on the first 3 principal components as discussed in section 4. The second method is the random walk. The third method is the mean forecast, when the term structure a year ahead is predicted to be equal to the average term structure so far. Finally, we consider Diebold-Li forecasting procedure.

The Diebold and Li's (2002) procedure consists of the following steps. First, we regress the term structure on three deterministic curves, the components of the Nelson and Siegel's (1987) forward rate curve:

$$f_t(T) = \beta_{1t} + \beta_{2t} e^{-\lambda_t T} + \beta_{3t} \lambda_t T e^{-\lambda_t T}.$$

This regression is run for each day in the subsample which the forecast is based upon.[4] Then, the time series for the coefficients of the regression are modeled as 3 separate autoregressive processes of order 1 (each of the the current coefficient is regressed on the corresponding coefficient one year before). A one-year ahead forecast of the coefficients is made, and the corresponding Nelson-Siegel forward curve is taken as the one-year ahead forecast of the term structure.

Figure 4 shows the predictive performance of the alternative methods considered.

---

[4] We fix parameter $\lambda_t$ at 0.0609 as Diebold and Li (2002) do.

**Figure 4**. The Predictive Performance Of Different Forecasting Methods

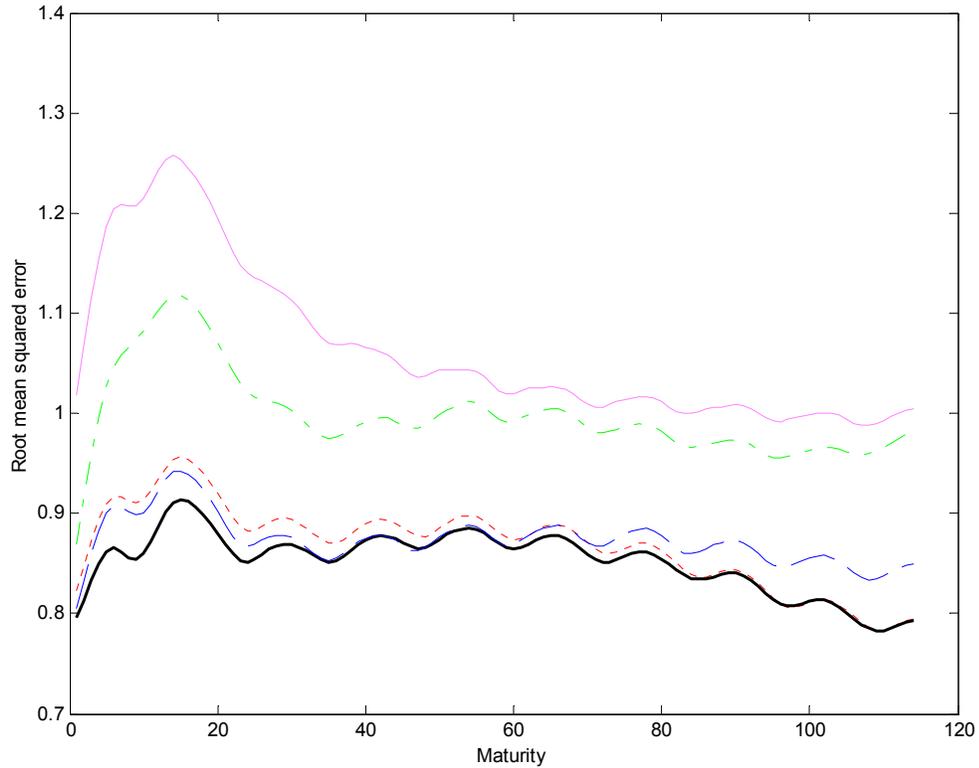

The thick solid line on the above graph corresponds to our method. The dashed line is for functional autoregression estimated with principal components. The dotted line is for mean prediction. The dash-dot line is for Diebold and Li (2002) method and solid thin line is for random walk. Our method has the best pseudo out of sample forecasting record uniformly across different maturities. The functional autoregression estimated with principal components is the second best for relatively short maturities and the third best for long maturities where the mean prediction works equally well with our preferred method. The worst performance is shown by the random walk.

It should be stressed that this ranking of the forecasting performances crucially depends on whether we exclude the recession period from the consideration. If we include the recession period, random walk outperforms all the other methods.

**Conclusion**

**Appendix A**

Consider an abstract real Hilbert space $H$. Let function $f_n$ map a probability space $(\Omega, \mathbb{A}, P)$ to $H$. We call this function an H-valued random variable if the scalar product $(g, f_n)$ is a standard random variable for every $g$ from $H$.[5]

**Definition 1**. If $E\|f\| < \infty$, then there exists an element of $H$, denoted as $Ef$ and called *expectation* of $f$, such that
$$E(g, f) = (g, Ef), \text{ for any } g \in H.$$

**Definition 2**. Let $f$ be an H-valued random variable, such that $E\|f\|^2 < \infty$ and $Ef = 0$. The *covariance operator* of $f$ is the bounded linear operator on $H$, defined by
$$C_f(g) = E[(g, f)f], \quad g \in H.$$

If $Ef \neq 0$, one sets $C_f = C_{f-Ef}$.

**Definition 3.** Let $f_1$ and $f_2$ be two H-valued random variables, such that $E\|f_1\|^2 < \infty, E\|f_2\|^2 < \infty$ and $Ef_1 = Ef_2 = 0$. Then the *cross-covariance operators* of $f_1$ and $f_2$ are bounded linear operators on $H$ defined by
$$C_{f_1, f_2}(g) = E[(g, f_1)f_2], \quad g \in H,$$
$$C_{f_2, f_1}(g) = E[(g, f_2)f_1], \quad g \in H.$$

If $Ef_1 \neq 0$ or/and $Ef_2 \neq 0$, one sets
$$C_{f_1, f_2} = C_{f_1 - Ef_1, f_2 - Ef_2},$$
$$C_{f_2, f_1} = C_{f_2 - Ef_2, f_1 - Ef_1}.$$

**Definition 4.** A sequence $\{\eta_n, \ n \in Z\}$ of H-valued random variables is said to be *H-white noise* if
1) $0 < E\|\eta_n\|^2 = \sigma^2 < \infty, \quad E\eta_n = 0, \quad C_{\eta_n}$ do not depend on $n$ and
2) $\eta_n$ is orthogonal to $\eta_m$; $n, m \in Z, \quad n \neq m$; i.e.,
$$E\{(x, \eta_n)(y, \eta_m)\} = 0, \text{ for any } x, y \in H.$$
$\{\eta_n, \ n \in Z\}$ is said to be a **strong H-white noise** if it satisfies 1) and
2') $\{\eta_n, \ n \in Z\}$ is a sequence of i.i.d. H-valued random variables.

**Example: Stochastic Processes**

Consider a set of stochastic processes, $f_i(T)$, on interval $[0, \overline{T}]$. Let $Ef_i(T) = 0$ for each $T$. Let the covariance function of each process be $Ef_i(S)f_i(T) = \Gamma_{ii}(S, T)$, and cross-covariance function between two processes be $Ef_i(S)f_j(T) = \Gamma_{ij}(S, T)$. Assume that with probability 1 the sample paths of the processes are in $L^2[0, \overline{T}]$. Each stochastic process defines an H-valued random variable with zero mean. The covariance operator of $f_i$ is the integral operator with kernel $\Gamma_{ii}(S, T)$:
$$C_{f_i}: g(T) \mapsto \int_0^{\overline{T}} \Gamma_{ii}(S, T) g(T) dT,$$
and the cross-covariance operator of $f_i$ and $f_j$ is the integral operator with kernel $\Gamma_{ij}(S, T)$.

---

[5] The definitions that follow are slight modifications of those in Chapters 2,3 of Bosq (2000).

**Appendix B**
**Proof of Theorem 2:** Transform the objective function in problem (RR) as

$$E\|f_{t+1} - AB'f_t\|^2 = tr(\Gamma_{11} - AB'\Gamma_{21} - \Gamma_{12}BA' + AB'\Gamma_{11}BA')$$
$$= tr(\Gamma_{11}) - tr(AB'\Gamma_{21} + \Gamma_{12}BA') + tr(AA')$$
$$= tr(\Gamma_{11}) - 2tr(B'\Gamma_{21}A) + tr(A'A),$$

where the first equality follows from the fact that the expectation of the squared norm of an $L^2$-valued random variable is equal to the trace of its covariance operator (see Bosq (2000) p.37), and the second equality follows from the constraint $B'\Gamma_{11}B = I_k$ imposed on $B$. To see that the third equality holds, write $tr(AA') = \sum_{i=1}^{\infty} e_i'AA'e_i$ and $tr(AB'\Gamma_{21}) = \sum_{i=1}^{\infty} e_i'AB'\Gamma_{21}e_i$, where $\{e_i\}$ is an arbitrary basis in $L^2$. Then use the fact that $A$ and $B'\Gamma_{21}$ are finite-dimensional vectors of functions from $L^2$, and apply Parceval's equality.

We will first minimize the transformed objective function with respect to $A$, taking $B$ as given. Note, that since $A = [a_1 \ldots a_k]$ and $B = [b_1 \ldots b_k]$, the objective function can be further transformed to the form $tr(\Gamma_{11}) - 2\sum_{i=1}^{k} b_i'\Gamma_{21}a_i + \sum_{i=1}^{k} a_i'a_i$. Hence, the problem separates into $k$ different problems: $-2b_i'\Gamma_{21}a_i + a_i'a_i \to \min$.

A necessary condition for the optimal $a_i$ to exist is that the Fréchet derivative of the objective function with respect to $a_i$ is equal to zero (see, for example, proposition 2 in §7.2 and theorem 1 in §7.4 of Luenberger (1969)). That is, $-2\Gamma_{12}b_i + 2a_i = 0$ and we have $A = \Gamma_{12}B$, which proves statement ii) of theorem 2.

Substituting $A = \Gamma_{12}B$ into the objective function, we reformulate problem (RR) as $\sum_{i=1}^{k} b_i'\Gamma_{21}\Gamma_{12}b_i \to \max$, subject to a constraint $B'\Gamma_{11}B = I_k$ and a requirement that $B'\Gamma_{21}\Gamma_{12}B$ is a diagonal $k \times k$ matrix with decreasing elements along the diagonal. This is equivalent to the recursive problem, formulated in statement iii) of the theorem. The mean square error of prediction that uses $k$ predictive factors is equal to $tr(\Gamma_{11}) - \max \sum_{i=1}^{k} b_i'\Gamma_{21}\Gamma_{12}b_i$. Hence, the $\max b_i'\Gamma_{21}\Gamma_{12}b_i$ can be interpreted as the reduction in the mean square error of forecasting due to the i-th predictive factor as iv) claims.

Now, let us define $x_i = \Gamma_{11}^{1/2}b_i$. Since $\Gamma_{12} = \rho\Gamma_{11}$, we have

$$b_i'\Gamma_{21}\Gamma_{12}b_i = b_i'\Gamma_{11}\rho'\rho\Gamma_{11}b_i$$
$$= x_i'\Gamma_{11}^{1/2}\rho'\rho\Gamma_{11}^{1/2}x_i,$$

and the recursive problem can be reformulated in terms of $x_i$ as

$$x_i'\Gamma_{11}^{1/2}\rho'\rho\Gamma_{11}^{1/2}x_i \to \max \text{ s.t. } x_i'x_j = 0, j < i \text{ and } x_i'x_i = 1.$$

Since $\Gamma_{11}^{1/2}\rho'\rho\Gamma_{11}^{1/2}$ is a compact[6] self-adjoint operator, a solution $x_i$ to this problem exists and equal to the eigenvector corresponding to the $i$-th largest eigenvalue $\lambda_i$ of $\Gamma_{11}^{1/2}\rho'\rho\Gamma_{11}^{1/2}$. The maximum value of the objective function is equal to $\lambda_i$ (see the proof of the spectral theorem III.5.1 in Gohberg and Gohberg (1981)).

---

[6] Compactness follows from the fact that $\Gamma_{11}^{1/2}$ is compact and the product of a compact operator and a bounded operator is compact.

To prove *i)* and *iv)*, it remains to note that $\Gamma_{11}^{-1/2} x_i$ are eigenvectors of $\Gamma_{21}\Gamma_{12} - \lambda \Gamma_{11}$ corresponding to the eigenvalues $\lambda_i$ and that any eigenvector of the above operator pencil has form $\Gamma_{11}^{-1/2} x_i$.[7] Indeed,

$$(\Gamma_{21}\Gamma_{12} - \lambda_i \Gamma_{11})\Gamma_{11}^{-1/2} x_i = (\Gamma_{11}^{1/2}\Gamma_{11}^{1/2} \rho' \rho \Gamma_{11}^{1/2} - \lambda_i \Gamma_{11}^{1/2}) x_i$$
$$= \Gamma_{11}^{1/2}(\Gamma_{11}^{1/2} \rho' \rho \Gamma_{11}^{1/2} - \lambda_i I) x_i$$
$$= 0.$$

Therefore, $\Gamma_{11}^{-1/2} x_i$ is the eigenvector. On the other hand, suppose $y$ is an eigenvector of $\Gamma_{21}\Gamma_{12} - \lambda \Gamma_{11}$, then

$$(\Gamma_{21}\Gamma_{12} - \lambda \Gamma_{11})y = \Gamma_{11}^{1/2}(\Gamma_{11}^{1/2} \rho' \rho \Gamma_{11}^{1/2} - \lambda I)\Gamma_{11}^{1/2} y = 0,$$

which can be true only if $\Gamma_{11}^{1/2} y$ is an eigenvector of $\Gamma_{11}^{1/2} \rho' \rho \Gamma_{11}^{1/2}$ because $\Gamma_{11}^{1/2}$ is positive definite.

Under conditions of *v)*, the eigenvectors $x_i$ corresponding to the first $k$ eigenvalues of $\Gamma_{11}^{1/2} \rho' \rho \Gamma_{11}^{1/2}$ are uniquely determined (up to the sign) and $\Gamma_{11}^{1/2}$ has a trivial null space, so that $\Gamma_{11}^{-1/2}$ is uniquely defined on the range of $\Gamma_{11}^{1/2}$. Hence, the first $k$ eigenvectors of $\Gamma_{21}\Gamma_{12} - \lambda \Gamma_{11}$ are uniquely determined (up to a change in the sign) and *v)* follows.□

**Appendix C**

**Proof of Theorem 3:** Consider Rayleigh functionals:

$$\gamma(b) = \frac{b'\Gamma_{21}\Gamma_{12}b}{b'\Gamma_{11}b}, \quad \gamma_\alpha(b) = \frac{b'\Gamma_{21}\Gamma_{12}b}{b'(\Gamma_{11} + \alpha I)b}, \text{ and } \hat{\gamma}_\alpha(b) = \frac{b'\hat{\Gamma}_{21}\hat{\Gamma}_{12}b}{b'(\hat{\Gamma}_{11} + \alpha I)b}$$

for operator pencils $\Gamma_{21}\Gamma_{12} - \lambda \Gamma_{11}$, $\Gamma_{21}\Gamma_{12} - \lambda(\Gamma_{11} + \alpha I)$, and $\hat{\Gamma}_{21}\hat{\Gamma}_{12} - \lambda(\hat{\Gamma}_{11} + \alpha I)$ respectively. According to the max-min principle (see Eschwé and Langer (2004)), the eigenvalues of the above operator pencils satisfy relationships:

$$\lambda_j = \max_{\dim M = j} \min_{b \in M} \gamma(b), \quad \lambda_{\alpha j} = \max_{\dim M = j} \min_{b \in M} \gamma_\alpha(b), \text{ and } \hat{\lambda}_{\alpha j} = \max_{\dim M = j} \min_{b \in M} \hat{\gamma}_\alpha(b).$$

To prove convergence $\hat{\lambda}_{\alpha j} \to \lambda_j$, we will explore relationship between these Rayleigh functionals. We have the following

**Proposition 1.** *Suppose that* $\alpha \to 0$, $\left(\dfrac{n}{\log n}\right)^{1/2} \alpha \to \infty$ *as* $n \to \infty$. *Then*

$$\sup_{b \in L^2} |\hat{\gamma}_\alpha(b) - \gamma_\alpha(b)| \to 0 \text{ in probability as } n \to \infty.$$

The proof is based on an extension of Lemma 1 in Leurgans et al (1993) that we formulate in Lemma 1 below. Given the extension, the proof is essentially the same as that of Proposition 3 in Leurgans et al (1993) and we omit it here.

To formulate Lemma 1, let us define $\Delta_1^{(n)} = \hat{\Gamma}_{11} - \Gamma_{11}$ and $\Delta_2^{(n)} = \hat{\Gamma}_{21}\hat{\Gamma}_{12} - \Gamma_{21}\Gamma_{12}$. We define the norm of integral operator $A$ as $\|A\|^2 = \iint A(s,t)^2 \, ds \, dt$. That is to say, $\|A\|$ is the Hilbert-Schmidt norm of $A$. We will use the property of the Hilbert-Schmidt norm that $\|AB\| \leq \|A\| \cdot \|B\|$. The following Lemma is a key to the proof of Proposition 1:

---

[7] Note that although $\Gamma_{11}^{-1/2}$ is unbounded, it can be defined on $\text{Im}(\Gamma_{11}^{1/2})$, which includes the eigenvectors of $\Gamma_{11}^{1/2} \rho' \rho \Gamma_{11}^{1/2}$.

**Lemma 1.** Let $\delta_n = \max(\|\Delta_1^{(n)}\|, \|\Delta_2^{(n)}\|)$. Then under Assumption 2, $\delta_n = O_P\left(\left(\frac{\log n}{n}\right)^{1/2}\right)$.

**Proof:** Theorems 4.1 and 4.8 of Bosq (2000) imply that $\|\Delta_1^{(n)}\| = O_P\left(\left(\frac{\log n}{n}\right)^{1/2}\right)$ and that

$\|\hat{\Gamma}_{12} - \Gamma_{12}\|$ and $\|\hat{\Gamma}_{21} - \Gamma_{21}\|$ are $O_P\left(\left(\frac{\log n}{n}\right)^{1/2}\right)$. Now,

$$\|\hat{\Gamma}_{21}\hat{\Gamma}_{12} - \Gamma_{21}\Gamma_{12}\| = \|\hat{\Gamma}_{21}\hat{\Gamma}_{12} - \Gamma_{21}\hat{\Gamma}_{12} + \Gamma_{21}\hat{\Gamma}_{12} - \Gamma_{21}\Gamma_{12}\|$$

$$\leq \|\hat{\Gamma}_{21}\hat{\Gamma}_{12} - \Gamma_{21}\hat{\Gamma}_{12}\| + \|\Gamma_{21}\hat{\Gamma}_{12} - \Gamma_{21}\Gamma_{12}\|$$

$$\leq \|\hat{\Gamma}_{21} - \Gamma_{21}\|\|\hat{\Gamma}_{12}\| + \|\Gamma_{21}\|\|\hat{\Gamma}_{12} - \Gamma_{12}\|$$

$$= O_P\left(\left(\frac{\log n}{n}\right)^{1/2}\right),$$

which completes the proof.□

Using Proposition 1, it is easy to prove part i) of our theorem. Note that since $\gamma(b) \geq \gamma_\alpha(b)$ for any $b \in L^2$, we have:

$$\lambda_j = \max_{\dim M = j} \min_{b \in M} \gamma(b) \geq \lambda_{\alpha j} = \max_{\dim M = j} \min_{b \in M} \gamma_\alpha(b) \geq \min_{b \in sp(b_1,\ldots,b_j)} \gamma_\alpha(b).$$

But $\gamma_\alpha(b) \xrightarrow{\alpha \to 0} \gamma(b)$ uniformly over any finite dimensional subspace $M \subset L^2$. This is because 1) the ratio $\gamma_\alpha(b)/\gamma(b)$ differs from 1 by $(\alpha b'b)/(b'\Gamma_{11}b)$, 2) $\Gamma_{11}$ is positive definite and therefore the minimum of $(b'\Gamma_{11}b)/(b'b)$ over ant finite-dimensional $M$ is strictly greater than zero, and 3) $\max_{b \in M}|(\alpha b'b)/(b'\Gamma_{11}b)| \to 0$ as $\alpha \to 0$.

Hence

$$\lambda_j \geq \lambda_{\alpha j} \geq \min_{b \in sp(b_1,\ldots,b_j)} \gamma_\alpha(b) \xrightarrow{\alpha \to 0} \min_{b \in sp(b_1,\ldots,b_j)} \gamma(b) = \gamma(b_j) = \lambda_j,$$

which implies that $\lambda_{\alpha j} \xrightarrow{\alpha \to 0} \lambda_j$.

Furthermore, Proposition 1 implies that

$$\left|\lambda_{\alpha j} - \hat{\lambda}_{\alpha j}\right| = \left|\max_{\dim M = j} \min_{b \in M} \gamma_\alpha(b) - \max_{\dim M = j} \min_{b \in M} \hat{\gamma}_\alpha(b)\right| \xrightarrow{p} 0.$$

Thus,

$$\left|\lambda_j - \hat{\lambda}_{\alpha j}\right| \leq \left|\lambda_j - \lambda_{\alpha j}\right| + \left|\lambda_{\alpha j} - \hat{\lambda}_{\alpha j}\right| \xrightarrow{p} 0,$$

which proves part i) of the theorem.

Let us now turn to part ii) of the theorem. In our proof we will use the following fact:

(A1) $$\lambda_j^{-1} \gamma_\alpha(\hat{b}_{\alpha j}) \xrightarrow{p} 1.$$

It follows from part i) of the theorem, from Proposition 1, and from the inequality

$$\left|\gamma_\alpha(\hat{b}_{\alpha j}) - \lambda_j\right| = \left|\gamma_\alpha(\hat{b}_{\alpha j}) - \hat{\gamma}_\alpha(\hat{b}_{\alpha j}) + \hat{\lambda}_{\alpha j} - \lambda_j\right| \leq \left|\gamma_\alpha(\hat{b}_{\alpha j}) - \hat{\gamma}_\alpha(\hat{b}_{\alpha j})\right| + \left|\hat{\lambda}_{\alpha j} - \lambda_j\right|.$$

Let us define $\beta_{ji} = \hat{b}_{\alpha j}'\Gamma_{11}b_i$, $i \le j$. These can be interpreted as the coefficients in the least squares regression of $\hat{b}_{\alpha j}$ on $b_{\alpha 1},\ldots,b_{\alpha j}$ in metric $\Gamma_{11}$. Define the error of such a regression as $s_i$, that is $\hat{b}_{\alpha j} = \sum_{i=1}^{j} \beta_{ji} b_i + s_j$. Note that $s_j'\Gamma_{11}b_i = 0$ and, since $b_i$ are eigenvectors of $\Gamma_{21}\Gamma_{12} - \lambda\Gamma_{11}$, $s_j'\Gamma_{21}\Gamma_{12}b_i = 0$ too.

Consider first the case $j = 1$. We have:
$$\lambda_1^{-1}\gamma_\alpha(\hat{b}_{\alpha 1}) = \frac{\beta_{11}^2 + \lambda_1^{-1}s_1'\Gamma_{21}\Gamma_{12}s_1}{\beta_{11}^2 + s_1'\Gamma_{11}s_1 + \alpha\hat{b}_{\alpha 1}'\hat{b}_{\alpha 1}}.$$

Eigenvalues of the operator pencil $\Gamma_{21}\Gamma_{12} - \lambda\Gamma_{11}$ must satisfy $\lambda_j = \max_{b \perp_{\Gamma_{11}} sp(b_1,\ldots,b_{j-1})} \frac{b'\Gamma_{21}\Gamma_{12}b}{b'\Gamma_{11}b}$.

Therefore, $\frac{s_1'\Gamma_{21}\Gamma_{12}s_1}{s_1'\Gamma_{11}s_1} \le \lambda_2$ and we have:
$$\lambda_1^{-1}\gamma_\alpha(\hat{b}_{\alpha 1}) \le \frac{1 + (\lambda_2/\lambda_1)\beta_{11}^{-2}s_1'\Gamma_{11}s_1}{1 + \beta_{11}^{-2}s_1'\Gamma_{11}s_1 + \alpha\beta_{11}^{-2}\hat{b}_{\alpha 1}'\hat{b}_{\alpha 1}}.$$

The left hand side of the above inequality tends to 1 by (A1), and the right hand side is bounded from above by 1. Hence the right hand side must tend to 1 and we have:

(A2) $\qquad \beta_{11}^{-2}s_1'\Gamma_{11}s_1 \xrightarrow{p} 0$ and $\alpha\beta_{11}^{-2}\hat{b}_{\alpha 1}'\hat{b}_{\alpha 1} \xrightarrow{p} 0$

Note that the first part of (A2), definition of $s_1$, and the normalization $b_1'\Gamma_{11}b_1 = 1$ imply that $\beta_{11}^{-2}\hat{b}_{\alpha 1}'\Gamma_{11}\hat{b}_{\alpha 1} \xrightarrow{p} 1$. Recalling that $\hat{b}_{\alpha 1}$ is normalized so that $\hat{b}_{\alpha 1}'(\hat{\Gamma}_{11} + \alpha I)\hat{b}_{\alpha 1} = 1$, (A2) therefore implies that
$$\beta_{11}^{-2}\hat{b}_{\alpha 1}'(\hat{\Gamma}_{11} - \Gamma_{11})\hat{b}_{\alpha 1} + \beta_{11}^{-2} \xrightarrow{p} 1.$$

But
$$\beta_{11}^{-2}\hat{b}_{\alpha 1}'(\hat{\Gamma}_{11} - \Gamma_{11})\hat{b}_{\alpha 1} \le \alpha\beta_{11}^{-2}\hat{b}_{\alpha 1}'\hat{b}_{\alpha 1}\frac{\delta_n}{\alpha} \to 0$$

by Lemma 1 and the assumption of the theorem that $\alpha(n/\log n)^{1/2} \to \infty$. Therefore, $\beta_{11}^2 \xrightarrow{p} 1$. Since we would like to prove all convergences up to a choice of the sign of the vectors in the converging series, we will ignore the possibility of $\beta_{11} \xrightarrow{p} (-1)$ and will only consider the case $\beta_{11} \xrightarrow{p} 1$. This fact together with (A2) imply that
$$(\hat{b}_{\alpha j} - b_j)'\Gamma_{11}(\hat{b}_{\alpha j} - b_j) \xrightarrow{p} 0 \text{ and } \alpha\hat{b}_{\alpha 1}'\hat{b}_{\alpha 1} \xrightarrow{p} 0.$$

Suppose now that statement ii) of the theorem is true for $j \le k-1$. Suppose also that for $j \le k-1$, $\alpha\hat{b}_{\alpha j}'\hat{b}_{\alpha j} \xrightarrow{p} 0$. We will prove that $\alpha\hat{b}_{\alpha k}'\hat{b}_{\alpha k} \xrightarrow{p} 0$ and that statement ii) of the theorem holds for $j = k$.

First, we will show that $\beta_{kj} \xrightarrow{p} 0$, $j \le k-1$. Because of the normalization $\hat{b}_{\alpha k}'(\hat{\Gamma}_{11} + \alpha I)\hat{b}_{\alpha j} = \delta_{kj}$, where $\delta_{kj} = 0$ if $k \ne j$ and $\delta_{kj} = 0$ if $k = j$, we have:
$$\beta_{kj} = \hat{b}_{\alpha k}'\Gamma_{11}b_j$$
$$= \hat{b}_{\alpha k}'\Gamma_{11}(b_j - \hat{b}_{\alpha j}) + \hat{b}_{\alpha k}'(\Gamma_{11} - \hat{\Gamma}_{11})\hat{b}_{\alpha j} - \alpha\hat{b}_{\alpha k}'\hat{b}_{\alpha j},$$

absolute value of which is less than or equal to

$$\sqrt{\hat{b}_{\alpha k}'\Gamma_{11}\hat{b}_{\alpha k}}\sqrt{(b_j - \hat{b}_{\alpha j})\Gamma_{11}(b_j - \hat{b}_{\alpha j})} + \sqrt{\alpha\hat{b}_{\alpha k}'\hat{b}_{\alpha k}}\sqrt{\alpha\hat{b}_{\alpha j}'\hat{b}_{\alpha j}}\left(\frac{\delta_n}{\alpha} + 1\right).$$

The second term in the above expression tends in probability to zero because $\alpha\hat{b}_{\alpha j}'\hat{b}_{\alpha j} \xrightarrow{p} 0$ by assumption of our induction, whereas $\alpha\hat{b}_{\alpha k}'\hat{b}_{\alpha k} < 1$ by the normalization assumption and $\delta_n/\alpha \to 0$ by assumption of the theorem. The first term also tends to zero. It is because $(b_j - \hat{b}_{\alpha j})\Gamma_{11}(b_j - \hat{b}_{\alpha j}) \xrightarrow{p} 0$ by assumption of the induction and $\hat{b}_{\alpha k}'\Gamma_{11}\hat{b}_{\alpha k}$ is bounded in probability. Indeed, $\hat{b}_{\alpha k}'\Gamma_{11}\hat{b}_{\alpha k} = \hat{b}_{\alpha k}'(\Gamma_{11} - \hat{\Gamma}_{11})\hat{b}_{\alpha k} + \hat{b}_{\alpha k}'\hat{\Gamma}_{11}\hat{b}_{\alpha k}$, where the first term in the right hand side tends in probability to zero, and the second term is less than 1.

Let us now show that $\beta_{kk}^2 \xrightarrow{p} 1$. For the $k+1$-th eigenvalue of $\Gamma_{21}\Gamma_{12} - \lambda\Gamma_{11}$, we have:

$$\lambda_{k+1} = \max_{b \perp_{\Gamma_{11}} sp(b_1,\ldots,b_k)} \frac{b'\Gamma_{21}\Gamma_{12}b}{b'\Gamma_{11}b} \geq \frac{s_k'\Gamma_{21}\Gamma_{12}s_k}{s_k'\Gamma_{11}s_k}.$$

Expanding the right hand side of the above inequality, we have:

$$\lambda_{k+1} \geq \frac{\hat{b}_{\alpha k}'\Gamma_{21}\Gamma_{12}\hat{b}_{\alpha k} - \lambda_1\beta_{k1}^2 - \ldots - \lambda_k\beta_{kk}^2}{\hat{b}_{\alpha k}'\Gamma_{11}\hat{b}_{\alpha k} - \beta_{k1}^2 - \ldots - \beta_{kk}^2}$$

$$\geq \frac{\hat{b}_{\alpha k}'\Gamma_{21}\Gamma_{12}\hat{b}_{\alpha k} - \lambda_1\beta_{k1}^2 - \ldots - \lambda_k\beta_{kk}^2}{\hat{b}_{\alpha k}'(\Gamma_{11} + \alpha I)\hat{b}_{\alpha k} - \beta_{k1}^2 - \ldots - \beta_{kk}^2},$$

Using Proposition 1, part i) of the theorem, and the just proven fact that $\beta_{kj} \xrightarrow{p} 0$, $j \leq k-1$, we have:

$$\lambda_{k+1} \geq \lambda_k \frac{1 - \beta_{kk}^2 + o_p(1)}{1 - \beta_{kk}^2 + o_p(1)}.$$

But $\lambda_{k+1} < \lambda_k$ by assumption 3. Hence there exists $\varphi < 1$, such that

$$\frac{1 - \beta_{kk}^2 + o_p(1)}{1 - \beta_{kk}^2 + o_p(1)} < \varphi$$

with probability 1. This only can be true if $\beta_{kk}^2 \xrightarrow{p} 1$.

Now, let us conclude the step of the induction. We have:

$$\lambda_k^{-1}\gamma_\alpha(\hat{b}_{\alpha k}) = \lambda_k^{-1}\frac{\lambda_1\beta_{k1}^2 + \ldots + \lambda_{k-1}\beta_{k,k-1}^2 + \lambda_k\beta_{kk}^2 + s_k'\Gamma_{21}\Gamma_{12}s_k}{\beta_{k1}^2 + \ldots + \beta_{k,k-1}^2 + \beta_{kk}^2 + s_k'\Gamma_{11}s_k + \alpha\hat{b}_{\alpha k}'\hat{b}_{\alpha k}}$$

$$\leq \frac{(\lambda_1/\lambda_k)\beta_{k1}^2 + \ldots + (\lambda_{k-1}/\lambda_k)\beta_{k,k-1}^2 + \beta_{kk}^2 + (\lambda_{k+1}/\lambda_k)s_k'\Gamma_{11}s_k}{\beta_{k1}^2 + \ldots + \beta_{k,k-1}^2 + \beta_{kk}^2 + s_k'\Gamma_{11}s_k + \alpha\hat{b}_{\alpha k}'\hat{b}_{\alpha k}}$$

$$= \frac{o_p(1) + 1 + (\lambda_{k+1}/\lambda_k)s_k'\Gamma_{11}s_k}{o_p(1) + 1 + s_k'\Gamma_{11}s_k + \alpha\hat{b}_{\alpha k}'\hat{b}_{\alpha k}}.$$

The left hand side of this inequality tends in probability to 1 by (A1). The right hand side is bounded in probability from above by $1+\varepsilon$ for any $\varepsilon > 0$. Hence, the right hand side tends in probability to 1 and $\alpha \hat{b}_{\alpha k}'\hat{b}_{\alpha k} \xrightarrow{p} 0$ and $s_k'\Gamma_{11}s_k \xrightarrow{p} 0$.

Since we would like to prove ii) only up to an appropriate sign of $\hat{b}_{\alpha k}$, we will assume that $\beta_{kk} \xrightarrow{p} 1$ and exclude the case $\beta_{kk} \xrightarrow{p} (-1)$ from our consideration. Then, $s_k'\Gamma_{11}s_k \xrightarrow{p} 0$ implies that $(\hat{b}_{\alpha k} - b_k)\Gamma_{11}(\hat{b}_{\alpha k} - b_k) \xrightarrow{p} 0$ and the induction is completed.□